\input amsppt1
\magnification=\magstep1
\documentstyle{amsppt}
\baselineskip 22pt
\pagewidth{13.5cm}
\pageheight{16.5cm}
\NoBlackBoxes
\TagsOnRight
\NoRunningHeads
\nologo
\define\ct{C_t}
\define\nx3{\nabla x_3}
\define\npx3{\nabla^\perp x_3}
\define\ve{\vert}
\define\ka{\kappa}
\define\la{\langle}
\define\ra{\rangle}
\define\pa{\partial}
\vskip 1.5cm
\topmatter
\title
On the total curvature of minimal annuli in R$^3$ 
and Nitsche's conjecture
\endtitle
\author
Qing CHEN
\endauthor
\affil
Department of Mathematics\\
The University of Science and technology of China\\
Hefei, Anhui, 230026\\
P.R.China\\
email: qchen\@sunlx06.nsc.ustc.edu.cn
\endaffil
\abstract 
{In this paper we prove the generalized Nitsche's conjecture
proposed by W. H. Meeks III and H. Rosenberg:
 For $t\ge0$, let $P_t$ denote the horizontal plane 
of height $t$ over the $x_1,x_2$-plane. Suppose that $M\subset R^3$ is a minimal
annulus with $\pa M\subset P_0$ and that $M$ intersects every $P_t$ in a 
simple closed curve. Then $M$ has finite total curvature. As a consequence,
we show that every properly embedded minimal surface of finite topology in $R^3$
with more than
one end has finite total curvature.}
\endabstract \endtopmatter

\subheading{ 1. Introduction}\par
The geometry and topology of minimal annuli in $R^3$ has attracted
geometer's attention for a long time.  For example,
B. Riemann described all minimal annuli in $R^3$ that are foliated by
the parallel circles, in terms of elliptic functions (see[9]); M. Shiffman
[10] studied minimal surfaces in $R^3$ bounded by 
two parallel convex curves.
In 1962, Nitsche [8] proved that a complete minimal annulus in $R^3$, which
intersects every horizontal plane in a star shaped curve, is a catenoid.
 Then he conjectured that every minimal embedded annulus in $R^3$ is
a catenoid, provided it intersects with every horizontal plane
in a simple closed curve.\par
In recent years W. H. Meeks III and his co-workers have  made much progress  
in understanding the structure of properly embedded minimal surface 
in $R^3$([3]-[7] etc.). In particular, It was shown in [5] and [6]
that Nitsche's conjecture is closely related to the finite total curvature
conjecture that every properly embedded minimal surface of finite topology 
 in $R^3$ with at least two ends
has finite total curvature (conjecture 1.1 of [6]). In [6], W. H. Meeks III
and H. Rosenberg reduced the finite total curvature conjecture to the
following generalized Nitsche's conjecture:
\proclaim{Conjecture } For $t\ge0$, let $P_t$ denote the horizontal plane 
of height $t$ over the $x_1,x_2$-plane. Suppose that $M\subset R^3$ is a minimal
annulus with $\pa M\subset P_0$ and that $M$ intersects every $P_t$ in a 
simple closed curve. Then $M$ has finite total curvature.\endproclaim\par
The aim of this paper is to prove the generalized Nitsche's conjecture.
Our proof is based on the curvature estimate and the careful analysis of geodesic
curvature and curvature of each $P_t\cap M$.  The fact that  each $P_t\cap M$ 
is an analytic curve is also fundamental in the proof.
\subheading {2. Curvature estimate}\par
Let $M$ be an immersed annulus in $R^3$ with the position vector 
$x=(x_1, x_2, x_3)$. The covariant derivatives of $R^3$ and $M$ 
are denoted by $D$ and $\nabla$ respectively. Let $P_t$ be the plane of 
$x_3=t$. 
From now on we suppose $\pa M\subset P_0$ and $M\cap P_t:=\ct$ is a single 
immersed curve, for every $t\in(0,+\infty)$. By Sard's theorem, for almost
of all $t\in(0,+\infty)$, $\nx3$ never vanish along $\ct$, and $\frac{\nx3}
{\ve\nx3\ve}$ is the normal vector field of $\ct$ in $M$.
\proclaim{Lemma 1}For almost of all $t>0$, $\ct$ is an analytic planner curve,
and its geodesic curvature $\ka_g$ and curvature $\ka$ are both real analytic.
Here we emphasize that $\ka$ is the curvature of $\ct$ as a plane curve.
\endproclaim
\demo{Proof} Suppose along $\ct$, $\nx3\ne0$. Giving a point $p\in\ct$, we 
choose an isothermal coordinate $(u,v)$ in a neighbourhood of $p$. Then
$x_3$ is a smooth harmonic and hence real analytic function of $(u,v)$. Without
lose of generality, suppose $\frac{\pa x_3}{\pa u}(p)\ne0$. Then the implicit
function theorem ([2] Chapter 10) reads that function $u=u(v)$ determined by
$x_3(u,v)=t$
is a real analytic function. Therefore locally $x\ve_{\ct}=x(v)$ is real
analytic; and the tangent vector of $\ct$, say $\frac{dx}
{dv}=x_u\frac{du}{dv}+x_v$, is analytic of $v$ and does not vanish. It follows
easily that $\ct$'s geodesic curvature and curvature are both real analytic with 
respect to $v$. Thus we complete the proof.  $\square$\enddemo
Suppose $T$ is the unit tangent vector field of $\ct$, such that $\{T,
\frac{\nx3}{\ve\nx3\ve}\}$ preserves the orientation of $M$. We have
$$\align
\ka_g&=\la D_TT,\frac{\nx3}{\ve\nx3\ve}\ra,\tag1\\
\ka_n&=\la D_TT,N\ra,\tag2
\endalign$$
where $N$ is the unit normal vector field of $M$ and $\ka_n$ is the normal
curvature of $\ct$. Putting $\npx3=Dx_3-\nx3$, the projection of constant
vector $Dx_3=e_3$ onto the normal direction of $M$, then (1) and (2) yield
$$\align \ve\nx3\ve\ka_g+\la\npx3,N\ra\ka_n&=\la D_TT,e_3\ra\\
&=-\la T,D_Te_3\ra\\
&=0.\tag3\endalign$$
Hence we obtain
$$\align
\ve\nx3\ve^2\ka_g^2&=\la\npx3,N\ra^2\ka_n^2\\
&=(1-\ve\nx3\ve^2)\ka_n^2.\tag4\endalign$$
Combinning the fact  $\ka^2=\ka_g^2+\ka_n^2$ with (4) we find 
$$\ka_g^2=(1-\ve\nx3\ve^2)\ka^2.\tag5$$
\par
Let $K$ denote the Gaussian curvature of $M$,  we have
\proclaim{Lemma 2} If $\nx3\ne0$ in $\ct$, then the following property holds:
$$\ve\ka\ve\le\frac{\sqrt{-K}}{\ve\nx3\ve}.\tag6$$
\endproclaim
\demo{Proof} Since $M$ is minimal immersed, $\sqrt{-K}$ is the absolute
value of the principal curvature of $M$. If $\ka_g=0$, then $\ve\ka\ve=
\ve\ka_n\ve
\le\sqrt{-K}$ and the lemma is valid. If $\ka_g\ne0$, we see from (4) and (5)
that
$$\align\ve\ka\ve&=\frac{\ve\ka_g\ve}{\sqrt{1-\ve\nx3\ve^2}}\\
&=\frac{\ve\ka_n\ve}{\ve\nx3\ve}\\
&\le\frac{\sqrt{-K}}{\ve\nx3\ve},\tag7\endalign$$
which completes the proof.   $\square$\enddemo
Although  $\ka$ and  $\ka_g$ might change the sign, by the analyticity 
we can derive from (5) that
\proclaim{Lemma 3} For almost of all $t$, along $\ct$ either $\ka_g=
\ka\sqrt{1-\ve\nx3\ve^2}$ or $\ka_g=-\ka\sqrt{1-\ve\nx3\ve^2}$ holds.
\endproclaim
\demo{Proof} From Lemma 1, $\ka$ and $\ka_g$ are both real analytic,
Lemma 3 is a direct consequence of (5) and an elementary lemma given in the
appendix.   $\square$\enddemo
\subheading {3. The proof of the generalized Nitsche's conjecture}\par
In this section we wish to prove the following theorem:
\proclaim{Theorem 1} Let $M$ be an embedded minimal annulus in $R^3$ such
that $\pa M\subset P_0$ and that $M$ intersects every $P_t(t>0)$ in a simple
closed curve. Then $M$ has finite total curvature.\endproclaim
To show the theorem we need
\proclaim{Lemma 4} $\int_{\ct}\ve\nx3\ve$ is independent of $t$.\endproclaim
\demo{Proof}Denote $M_t=M\cap\{0<x_3<t\}$. Since $x_3$ is a harmonic function,
Green's formula yields
$$\align 0&=\int_{M_t}\Delta x_3\\
&=\int_{\ct}\la\nx3,\frac{\nx3}{\ve\nx3\ve}\ra\quad-\int_{C_0}
                    \la\nx3,\frac{\nx3}{\ve\nx3\ve}\ra\\
&=\int_{\ct}\ve\nx3\ve -\int_{C_0}\ve\nx3\ve,  \endalign $$
this proves the lemma.   $\square$   \enddemo
\demo{Proof of Theorem 1} Put $R(t)=\int_{M_t}(-K)$, where $K$ is the 
Gaussian curvature of $M$. We have by the Gauss-Bonnet theorem
$$R(t)=\int_{\ct}\ka_g\quad+\int_{C_0}\ka_g.\tag8$$
Put $c_0=\int_{C_0}\ka_g$. Lemma 3 and (8) lead that for almost of
all $t>0$,
$$\align
R(t)&=c_0\pm\int_{\ct}\ka\sqrt{1-\ve\nx3\ve^2}\\
&=c_0\pm\int_{\ct}(\sqrt{1-\ve\nx3\ve^2}-1)\ka\pm\int_{\ct}\ka\\
&\le c_0+2\pi+\int_{\ct}(1-\sqrt{1-\ve\nx3\ve^2})\ve\ka\ve.\tag9
\endalign$$
Here the fact that $\int_{\ct}\ka=\pm2\pi$ is used in the last
inequality. Then (9) and Lemma 2 imply
$$\aligned
R(t)&\le c_1+\int_{\ct}\frac{\ve\nx3\ve^2}{1+\sqrt{1-\ve\nx3\ve^2}}\ve\ka\ve\\
&\le c_1+\int_{\ct}\ve\nx3\ve\sqrt{-K}\\
&\le c_1+\int_{\ct}\sqrt{-K},
\endaligned\tag10$$
where $c_1=c_0+2\pi$. By the Schwarz inequality, we derive from (10) 
that
$$\align R(t)&\le c_1+(\int_{\ct}\ve\nx3\ve)^{\frac12}(\int_{\ct}
           \frac{-K}{\ve\nx3\ve})^{\frac12}\\
&=c_1+c_2(\int_{\ct}\frac{-K}{\ve\nx3\ve})^{\frac12},\tag11
\endalign $$
where $c_2$ is another constant because of Lemma 4. 
The co-area formula([11]) leads
$R'(t)=\int_{\ct}\frac{-K}{\ve\nx3\ve}$, therefore (11) implies 
$$R(t)\le c_1+c_2\sqrt{R'(t)}\tag12$$
for almost of all $t\in(0,+\infty)$. \par
Because $R(t)$ is monotone non-decreasing of $t$, if there exists a $t_0>0$
such that $R(t)>c_1$ when $t\ge t_0$, by (12) we have for almost of all
$t>t_0$,
$$(R(t)-c_1)^2\le c_2^2R'(t),$$
or equivalently
$$\frac1{c_2^2}\le\frac{R'(t)}{(R(t)-c_1)^2}.\tag13$$
Integrating (13) from $t_0$ to $t$, one has
$$\align
\frac1{c_2^2}(t-t_0)&\le\frac1{R(t_0)-c_1}-\frac1{R(t)-c_1}\\
&\le\frac1{R(t_0)-c_1},  \endalign$$
which contradicts the fact that $t$ is unbounded. Thus we have proved the
theorem.   $\square$\enddemo
We notice that the embeddedness of each $\ct$ is only employed in estimate
(9). Actually we can show the following in the same way.
\proclaim{Theorem 2} Let $M$ be an immersed annulus in $R^3$. Suppose
for every $t>0$, $M\cap P_t$ is a closed plane curve with finite
 rotation index, and $\pa M\subset
P_0$. Then $M$ is of finite total curvature.\endproclaim

Theorem 1, together with Theorem 1.1 of [6], proves the finite total curvature
conjecture:
\proclaim{Theorem 3} Every properly embedded minimal surface of finite topology in $R^3$ with
more than one end has finite total curvature.\endproclaim
\subheading {4. Appendix} \par
\proclaim{Lemma 5} Let $f(t)$ and $g(t)$ be two real analytic functions in
$(-1,1)$, and $\lambda (t)$ a non-negative continuous function in $t$.
Suppose
$f^2=\lambda g^2$, then in $(-1, 1)$ either $f=\sqrt{\lambda}g$ or
$f=-\sqrt{\lambda}g$.\endproclaim
\demo{Proof} Suppose $f$ does not vanish identically, since the roots of $f$
is isolated, it suffices to prove
the lemma in a neighbourhood of a root of $f$. Without lose of 
generality, suppose $f(0)=0$ and $f\ne0$ otherwise in $(-\epsilon,\epsilon)$. We
express $f$ as
$$f(t)=t^nf_1(t),\qquad n>0\quad\text{and}\quad f_1(t)\ne0.$$ 
If $g(t)=t^mg_1(t)$ with $g_1(t)\ne0$, we see $m\le n$ and $\lambda(t)=
t^{2(n-m)}\frac{f_1^2(t)}{g_1^2(t)}$ is also real analytic
in $(-\epsilon,\epsilon)$.
Thus we can write $\lambda(t)=t^{2(n-m)}\lambda_1(t)$, $\lambda_1(t)\ne0$
in $(-\epsilon,\epsilon)$, and $f_1^2(t)=\lambda_1(t)g_1^2(t)$.
Therefore either $f_1(t)=\sqrt {\lambda_1(t)}g_1(t)$ or $f_1(t)=
-\sqrt {\lambda_1(t)}g_1(t)$ in $(-\epsilon,\epsilon)$. It follows
either $f=\sqrt{\lambda}g$ or $f=-\sqrt{\lambda}g$ in
$(-\epsilon,\epsilon)$. Which complete the proof of the assertion.  $\square$
\enddemo

\enddocument